\documentclass[a4paper,12pt]{article}

\usepackage{graphicx}
\usepackage{latexsym}
\usepackage{amsmath}
\usepackage{amsthm}
\usepackage{amsfonts}
\usepackage{amssymb}
\usepackage{enumerate}

\usepackage[authoryear]{natbib}
  4  \oddsidemargin-20pt \evensidemargin-20pt
\allowdisplaybreaks[4]

\newcommand\N{\mathbb{N}}

\newcommand\R{\mathbb{R}}

\newcommand{\vect}[1]{\mathbf{#1}}  %vectors
\newcommand\Prob{\mathbb{P}}    %probability
\newcommand\Exp{\mathbb{E}}     %expected value
\newcommand\ind{\mathbb{I}}     %indicators

\newcommand\Bb{\mathbb{B}}
\newcommand\Cb{\mathbb{C}}
\newcommand\Db{\mathbb{D}}
\newcommand\Gb{\mathbb{G}}
\newcommand\Hb{\mathbb{H}}
\newcommand\Sb{\mathbb{S}}
\newcommand\Tb{\mathbb{T}}

\DeclareMathOperator{\Cov}{Cov}

\DeclareMathOperator{\tKS}{KS}
\DeclareMathOperator{\tM}{M}
\DeclareMathOperator{\tL}{L}

\newcommand{\dceins}{ \dot C_1 }
\newcommand{\dczwei}{ \dot C_2 }
\newcommand{\dcp}{ \dot C_p }
\newcommand{\dchateins}{ \widehat{\dot C_1} }
\newcommand{\dchatzwei}{ \widehat{\dot C_2} }
\newcommand{\dchatp}{ \widehat{\dot C_p} }
\newcommand{\ec}{ C_n}

\newcommand\weak{\ \rightsquigarrow\ }
\newcommand\weakP{\ \overset{\Prob}{\underset{\xi}{\rightsquigarrow}}\ }

\newcommand\Pconv{ \stackrel{ \Prob}{\rightarrow} }

\newtheoremstyle{normal}% name
{2ex}               % Space above, empty = `usual value'
{3ex}               % Space below
{}                  % Body font
{}                  % Indent amount (empty = no indent, \parindent = para indent)
{\bfseries} % Thm head font
{}                  % Punctuation after thm head
{2pt}   % Space after thm head.
{\thmname{#1}\thmnumber{ #2.} \thmnote{(#3)}}% Thm head spec

\newtheoremstyle{italic}% name
{2ex}%      Space above, empty = `usual value'
{3ex}%      Space below
{\itshape}% Body font
{}%         Indent amount (empty = no indent, \parindent = para indent)
{\bfseries} % Thm head font
{}%        Punctuation after thm head
{2pt}%     Space after thm head.
{\thmname{#1}\thmnumber{ #2.} \thmnote{(#3)}}% Thm head spec

\theoremstyle{normal}
\newtheorem{definition}{Definition}[section]
\newtheorem{remark}[definition]{Remark}

\newtheorem{cond}[definition]{Condition}

\theoremstyle{italic}
\newtheorem{theorem}[definition]{Theorem}

\newtheorem{prop}[definition]{Proposition}

\begin{document}

\title{A test for Archimedeanity in bivariate copula models}

\author{Axel B\"ucher,	 Holger Dette  and  Stanislav Volgushev \\
Ruhr-Universit\"at Bochum \\
Fakult\"at f\"ur Mathematik \\
44780 Bochum, Germany \\
{\small e-mail: axel.buecher@ruhr-uni-bochum.de }\\
{\small e-mail: holger.dette@ruhr-uni-bochum.de}\\
{\small e-mail: stanislav.volgushev@ruhr-uni-bochum.de}\\
}

\maketitle

\begin{abstract}

We propose a new  test for the hypothesis that  a bivariate copula is an Archimedean copula.
 The test statistic is based on a combination of two measures resulting from the characterization of
Archimedean copulas by the property of associativity and by a strict upper bound on the diagonal by the
Fr\'{e}chet-upper bound.  We prove weak convergence of this statistic and show that
the critical values of the corresponding test can be  determined by the multiplier bootstrap
method. The test is shown to be consistent against all departures
from Archimedeanity if the copula satisfies weak smoothness assumptions.  A  simulation study is presented which
illustrates the  finite sample properties of the new test.

\end{abstract}


\section{Introduction}
\def\theequation{1.\arabic{equation}}
\setcounter{equation}{0}

Let $F$ be a bivariate continuous distribution function with marginal distribution
functions  $F_1$  and  $F_2$. By Sklar's Theroem [see \cite{sklar1959}] we
can decompose $F$ as follows
\begin{align}	\label{sklar}
	F(\vect{x}) = C(F_1(x_1),F_2(x_2)), \quad \vect{x}=(x_1,x_2)\in\R^2,
\end{align}
where $C$ is the unique copula associated to $F$. By definition, $C$ is a
bivariate distribution function on the unit square $[0,1]^2$ whose univariate
marginals are standard uniform distributions on the interval $[0,1]$. Equation
\eqref{sklar} is usually interpreted in the way that the copula $C$ completely
characterizes the information about the stochastic dependence contained in
$F$. For an extensive exposition on the theory of copulas we refer the reader
to the monograph \cite{nelsen2006}.

In the last decades, various parametric models for copulas have been developed,
among which the class of Archimedean copulas forms one the most famous and
largest class, see \cite{genemack1986, nelsen2006, mcnenesl2009} among many
others. Many widely used copulas, such as Clayton-, Gumbel- and Frank-copulas
are in fact Archimedean copulas. The elements of this class may be characterized
by a continuous, strictly decreasing  and convex function $\Phi:[0,1]\rightarrow[0,\infty]$
satisfying $\Phi(1)=0$ such that
\begin{align*}
	C(\vect{u})= \Phi^{[-1]}\left[\Phi(u_1)+\Phi(u_2)\right] \quad \text{for all } \vect{u}=(u_1,u_2)\in[0,1]^2.
\end{align*}
The function $\Phi$ is called the \textit{generator} of $C$ and its \textit{pseudo-ineverse}
$\Phi^{[-1]}(t)$ is defined as the usual inverse $\Phi^{-1}(t)$ for $t\in[0,\Phi(0)]$ and is
set to $0$ for $t\geq\Phi(0)$. The prominence of the class of Archimedean copulas
basically stems from the fact that they are  easy to handle and to simulate, see
\cite{genneszie2011}. While the estimation of Archimedean copulas has been
investigated in \cite{generive1993} and recently more thoroughly in \cite{genneszie2011},
the issue of testing for the hypothesis that the copula is an Archimedean one has found
much less interest in the literature. The present paper fills this gap by developing a
consistent test for this hypothesis.

Our interest in this problem stems from recent work of \cite{generive1993},
\cite{wangwell2000} and \cite{naifar2011} who proposed Archimedean copulas for modeling dependencies
between bivariate observations (among many others).
We also refer to the work of 
\cite{rivewell2001} who used
Archimedean copulas for modeling the dependence in the context of censored data.

To the best of our knowledge, the only available test hitherto has been discussed
in \cite{jaworski2010}. This author proposed a procedure which is based on a
characterization of Archimedean copulas similar to the one stated in Theorem 4.1.6
in \cite{nelsen2006} [which dates back to \cite{ling1965}]. To be precise recall that a
bivariate copula $C$ is called associative if and only if the identity
\begin{align} \label{asso}
	C(x,C(y,z))=C(C(x,y),z)\quad
\end{align}
holds for all $(x,y,z) \in [0,1]^3.$ Theorem 4.1.6 in \cite{nelsen2006} shows
that a bivariate copula  $C$  is an Archimedean copula if and only if $C$ is associative
and the inequality $C(u,u)<u$ holds for all $u\in(0,1)$, i.e. on the diagonal $C$ is strictly
dominated by the Fr\'{e}chet-upper bound $M(\vect{u})=\min(u_1,u_2)$. The procedure
suggested in \cite{jaworski2010} is in fact to test for associativity in order to check the
validity of an Archimedean copula model. The corresponding test statistic is defined as
\begin{align*}
	\mathcal{T}_n(x,y)=\sqrt{n}(\ec(x,\ec(y,y)) - \ec(\ec(x,y),y),
\end{align*}
where $(x,y)$ is some fixed point in the open cube $(0,1)^2$ and  $\ec$ denotes the
empirical copula, see Section \ref{sec:test} for details. The main advantage of this
approach is its simplicity, in particular the simple limit distribution of the resulting test
statistic, which is in fact normal. On the other hand this simplicity has its price
in terms of consistency. In  our opinion, the method proposed by \cite{jaworski2010}
has at least three mayor drawbacks. First of all, it is clearly not consistent against a
large class of alternatives since it only tests for equation \eqref{asso} with $y=z$.
Second,  \cite{jaworski2010} uses a pointwise approach in order to test for a global
hypothesis as in \eqref{asso}. This means that the test may not reject the hypothesis
because \eqref{asso} is satisfied at the particular point  $(x,y,y)$ under investigation,
although there may exist many other points where  \eqref{asso} is not satisfied. Third,
there exist copulas which are in fact associative but not Archimedean. These problems
also have strong implications for the practical applicability of the test as demonstrated
by results in a simulation study in \cite{jaworski2010}, where  the sample size has to
be chosen extremely large in order to get reasonable rejection probabilities.

To the best of our knowledge there exists no test for an Archimedean copula, which
is consistent against general alternatives and it is the primary purpose of this paper
to develop  such a procedure and to investigate its statistical properties. We propose
a test statistic which is based on a combination of two measures resulting from the
characterization of Archimedean copulas, namely the property of associativity as
described in \eqref{asso} and the  strict upper bound on the diagonal $C(u,u)<u$
for all $u\in(0,1)$. In Section \ref{sec:test}  we define a new process which is based
on an estimate of the difference of the left and right hand side of the defining equation
\eqref{asso} for associativity. We prove weak convergence of this process in the
space of all uniformly bounded functions on the cube $[0,1]^3$. As a consequence,
we also obtain weak  convergence of a corresponding  Cram\'{e}r-von-Mises and a Kolmogorov-Smirnov type statistic.
Because the asymptotic distribution depends in a complicated manner on the underlying copula we propose a multiplier bootstrap procedure to obtain the critical
values and show its validity. As a first main  result we obtain a test for associativity,
which is consistent against all alternatives satisfying weak smoothness assumptions
on $C$. In Section \ref{sec:testarch} we utilize these findings  to develop an asymptotic
test for the hypothesis of Archimedeanity. Finally in  Section \ref{sec:sim} we investigate
the finite sample performance of the new test by means of a simulation study.

\section{Testing Associativity}\label{sec:test}
\def\theequation{2.\arabic{equation}}
\setcounter{equation}{0}

\subsection{The test statistic and its asymptotic behavior}

In the following let $\vect{X}_1, \dots, \vect{X}_n$, $\vect{X}_i=(X_{i1},X_{i2})$
denote independent identically distributed bivariate random vectors with continuous distribution function $F$, marginal distribution functions $F_1$ and $F_2$
and copula $C=F(F_1^-, F_2^-)$. In this paragraph we will introduce a test statistic for the null hypothesis that the underlying copula
 is associative,
i.e. $C$ satisfies condition \eqref{asso} for all $(x,y,z) \in [0,1]^3$.

For this purpose we briefly  summarize relevant notations and results on the empirical copula,
which is the simplest and most popular  nonparametric estimator of  the copula. In particular
 we define
the empirical copula by
\begin{align*}
	\ec(\vect{u}) = F_n(F_{n1}^-(u_1), F_{n2}^-(u_2)),
\end{align*}
where $F_n (\vect{x}) =n^{-1}\sum_{i=1}^n \ind\{\vect{X}_i\leq\vect{x}\}$ and $F_{np } (x_p) =n^{-1}\sum_{i=1}^n \ind\{X_{ip}\leq x_p\}, p=1,2$ are the joint and marginal empirical distribution functions of the sample $\vect{X}_1, \dots, \vect{X}_n$, respectively. It is a well known result that
under the assumptions of  continuous partial derivatives of $C$ the corresponding empirical copula process
\begin{align}\label{ecp}
	\Cb_n = \sqrt{n}(\ec - C)
\end{align}
converges  weakly  towards a Gaussian limit field $\Gb_C$ in $l^\infty([0,1]^2)$, see \cite{ruschendorf1976, ferradweg2004, tsukahara2005} among others. Defining $\dcp$ as the $p$-th partial derivative of $C$ ($p=1,2$) the process $\Gb_C$ can be expressed as
\begin{align}\label{G_C}
	\Gb_C(\vect{x}) = \Bb_C(\vect{x}) - \dceins(\vect{x})\Bb_C(x_1,1) - \dczwei(\vect{x}) \Bb_C(1,x_2)
\end{align}
with the copula-brownian bridge $\Bb_C$, i.e. $\Bb_C$ is a centered Gaussian field with $\Cov(\Bb_C(\vect{x}),\Bb_C(\vect{y})) = C(\vect{x}\wedge\vect{y}) - C(\vect{x})C(\vect{y})$, where the minimum of two vectors is defined component-wise. As explained in \cite{segers2011} the assumption of continuity of the partial derivatives of $C$ on the whole unit square does not hold for many (even most) commonly used copula models and as a consequence Segers provides the result that the following nonrestrictive smoothness condition is sufficient
in order to obtain weak  convergence of the empirical copula process defined in \eqref{ecp}.

\begin{cond}\label{cond:pd}
\emph{For $p=1,2$ the first order partial derivative $\dcp$ of the copula $C$ with respect to $x_p$
exists and is continuous on the set $V_{p}=\{\vect{u}\in[0,1]^2\,:\, 0<u_p<1\}$.}
\end{cond}

Now, in order to test for associativity we consider the process
\begin{align*}
	\Hb_n(x,y,z)=\sqrt{n} \left\{\ec(x,\ec(y,z)) - \ec(\ec(x,y),z) \right\},
\end{align*}
where $(x,y,z)\in[0,1]^3$. The  asymptotic properties of the process $\{ \Hb_n(x,y,z)\}_{(x,y,z)\in [0,1]^3}$
are summarized in the following Theorem. Throughout this paper
$l^\infty( T) $ denotes the set of all uniformly bounded functions on $T$, and
the symbol $ \weak$ denotes
uniform convergence in a metric space (which will be specified in the corresponding statements).

\begin{theorem}\label{Hconv}
If the copula $C$ is associative and satisfies Condition \ref{cond:pd}, then it holds
$$\Hb_n \weak \Hb_C\quad\text{ in }~~ l^\infty([0,1]^3),$$
where the limit field $\Hb_C$ can be expressed as
\begin{align*}
	\Hb_C(x,y,z)=\Gb_C(x,C(y,z)) - \Gb_C(C(x,y),z) + \dczwei(x,C(y,z)) \Gb_C(y,z) - \dceins(C(x,y),z) \Gb(x,y).
\end{align*}
\end{theorem}

\textbf{Proof.}
If the copula $C$ is associative we can write the process $\Hb_n$ as
\begin{align*}
	\Hb_n = \sqrt{n} \left\{ \Phi(\ec) - \Phi(C) \right\},
\end{align*}
where the functional $\Phi:\Db_\Phi \rightarrow l^\infty([0,1]^3)$ is defined for
$$\alpha \in \Db_\Phi=\{F:F \text{ cdf on } [0,1]^2\} $$ by
\begin{align*}
	\Phi(\alpha)(x,y,z) = \alpha(x,\alpha(y,z)) - \alpha(\alpha(x,y),z).
\end{align*}
We will show later that under Condition \ref{cond:pd}  the mapping $\Phi$ is
Hadamard-differentiable at $C$ tangentially to the space
\begin{align*}
\Db_0=\left\{ \gamma \in C[0,1]^2\, | \,  \gamma(\vect{u})=0 \text{ for all } \, \vect{u} \in [0,1]^2 \text{ s.t. } C(\vect{u})\in\{0,1\}\right\},
\end{align*}
with derivative given by
\begin{align*}
	\Phi_C'(\alpha)(x,y,z) = \alpha(x,C(y,z)) - \alpha(C(x,y),z) + \dczwei(x,C(y,z)) \alpha(y,z) - \dceins(C(x,y),z) \alpha(x,y).
\end{align*}
Observing that $\Bb_C\in\Db_0$ a.s.,
the functional delta method, see Theorem 3.9.4 in \cite{vandwell1996}, yields the assertion.

We now briefly sketch how to see the Hadamard-differentiability of the mapping $\Phi$: let $t_n\rightarrow0$ and $\alpha_n\in l^\infty([0,1]^2)$ with $\alpha_n\rightarrow \alpha\in \Db_0$ such that $C+t_n\alpha_n\in\Db_\Phi$. Then
\begin{align*}
	t_n^{-1}\{ \Phi(C+t_n\alpha_n) - \Phi(C)\} = L_{n1} + L_{n2} - L_{n3}
\end{align*}
where
\begin{align*}
	L_{n1}(x,y,z) &= \alpha_n(x,(C+t_n\alpha_n)(y,z))  - \alpha_n((C+t_n\alpha_n)(x,y),z)  \\
	L_{n2}(x,y,z) &= t_n^{-1}\{ C(x,(C+t_n\alpha_n)(y,z)) - C(x,C(y,z))  \} \\
	L_{n3}(x,y,z) &= t_n^{-1}\{ C((C+t_n\alpha_n)(x,y),z) - C(C(x,y),z) \}.
\end{align*}
Exploiting the fact that $\alpha_n$ converges uniformly to a bounded function
and that $\alpha$ is uniformly continuous one can conclude that
$L_{n1}(x,y,z)= \alpha(x,C(y,z)) - \alpha(C(x,y),z) + o(1)$ uniformly in $(x,y,z)\in[0,1]^3$.
Regarding the summand $L_{n2}$ we have to split the investigation in two cases.
First, we consider all those $(x,y,z)\in[0,1]^3$ for which $C(y,z)\in(0,1)$.
A Taylor expansion of $C(x,\cdot)$ at $C(y,z)$ yields
$$L_{n2}(x,y,z)=\dczwei(x,C(y,z))\alpha_n(y,z) + r_n(x,y,z),$$
where the error term can be written as
$$r_n(x,y,z)=\big( \dczwei(x,u_n) - \dczwei(x,C(y,z)) \big)\alpha_n(y,z)$$
with some intermediate point $u_n$ between $C(y,z)$ and $(C+t_n\alpha_n)(y,z)$. The main term uniformly converges to $\dczwei(x,C(y,z)) \alpha(y,z)$ [note that partial derivatives of copulas are uniformly bounded by $1$] and it remains to show that $r_n (x,y,z)=o(1)$ uniformly in $(x,y,z)$ with $C(y,z)\in(0,1)$.

To see this,
we will show at the end of this proof  that for any $\varepsilon > 0$ there
exists a $\delta>0$, such that
\begin{align} \label{rest}
\limsup_{n\rightarrow\infty} \sup_{\vect{v} \in A_\delta} |\alpha_n(\vect{v})|\leq \varepsilon.
\end{align}
where $ \vect{v}=(y,z)$, $A_\delta =  \{ \vect{v}\in [0,1]^2 | \ C(\vect{v})\in[0,\delta)\cup(1-\delta,1]\}$. 
%Exploiting uniform convergence of $\alpha_n$,
%uniform continuity of $\alpha$ and the fact that $\alpha(\vect{v})=0$ for all
%$v\in A_0 = \{ \vect{v}\, | \, C(\vect{v})\in\{0,1\} \}$, we can conclude that there
%exists a $\delta>0$, such that
%$\limsup_{n\rightarrow\infty} \sup_{(y,z) \in A_\delta} |\alpha_n(y,z)|\leq \varepsilon$, where $A_\delta =  \{ \vect{v}\in [0,1]^2 | \ C(\vect{v})\in[0,\delta)\cup(1-\delta,1]\}$.
Then, since partial derivatives of copulas are bounded by $1$, we can conclude that
\begin{align*}
	\limsup_{n\rightarrow\infty} \sup_{x\in[0,1],(y,z) \in A_\delta} |r_n(x,y,z)| \leq \varepsilon.
\end{align*}
Due to Condition \ref{cond:pd} the partial derivative $\dczwei$ is uniformly
continuous on the quadrangle $[0,1]\times [\delta, 1-\delta]$. Thus, since
$\alpha$ is uniformly bounded and since $u_n\rightarrow C(y,z)$, we obtain
uniform convergence of $r_n(x,y,z)$ to $0$ for all $(y,z)$ s.t. $C(y,z)\in[\delta, 1-\delta]$,
i.e. for $(y,z)\in[0,1]^2\setminus A_\delta$. Combining the two facts derived above, it follows that
\begin{align*}
	\limsup_{n\rightarrow\infty} \sup_{x\in[0,1],C(y,z)\in(0,1)} |r_n(x,y,z)| \leq \varepsilon.
\end{align*}
Since $\varepsilon>0$ was arbitrary, this $\limsup$ must be zero. Summarizing, the case $(x,y,z)\in[0,1]^3$ such that $C(y,z)\in(0,1)$ is finished.

In the remaining case $C(y,z)\in\{0,1\}$, i.e. $(y,z)\in A_0$, Lipschitz-continuity of $C$ entails that
$$|L_{n2}(x,y,z)| = t_n^{-1}| C(x,C(y,z)+t_n\alpha_n(y,z)) - C(x,C(y,z)) | \leq \alpha_n(y,z) =\alpha(y,z)+o(1)=o(1)$$
uniformly in $(x,y,z)$ since in this case $\alpha(y,z)=0$. Finally, the summand $L_{n3}$ may be treated analogously.

To complete the proof it remains to show \eqref{rest}.
Exploiting uniform convergence of $\alpha_n$, uniform continuity of $\alpha$ and
the fact that $\alpha(\vect{v})=0$ for all $\vect{v}\in A_0 = \{ \vect{v}\, | \, C(\vect{v})\in\{0,1\} \}$,
we can conclude that there exists a $\kappa>0$ such that $|\alpha_n(\vect{v})| \leq \varepsilon$
for all $\vect{v}\in A_0^\kappa=\{\vect{v}\, | \, \exists \, \vect{u}\in A_0 \text{ s.t. } \|\vect{u}-\vect{v}\| \leq \kappa\}$ and sufficiently large $n$.
For $v_1\in[\kappa, 1]$ let $\delta(v_1)=\sup\{ C(v_1,z) \, | \, (v_1,z)\in A_0^\kappa \}$ [which equals $C(v_1,z(v_1))$
for some $z(v_1)$ such that $(v_1,z(v_1))\in  \overline{\partial A_0^\kappa\cap(0,1)^2}$ since for fixed any $v_1$ the function $u \mapsto C(v_1,u)$ i increasing]
and set $\delta=\inf_{v_1\in[\kappa,1]} \delta(v_1)$, which is strictly positive
due to compactness of  $\overline{\partial A_0^\kappa\cap(0,1)^2}$ and continuity of $C$. We will now show that this choice of $\delta$ yields \eqref{rest}.
Now, if $C(\vect{v})\leq \delta$, we have either $v_1< \kappa$ [then $\vect{v}\in A_0^\kappa$ since
$C(0,v_2)=0$] or $v_1\geq\kappa$. In the latter case, $C(\vect{v})\leq \delta(v_1)$ and monotonicity of $C$ imply
$\vect{v}\in A_0^\kappa$. This proves \eqref{rest} and completes the proof of Theorem \ref{Hconv}.
\qed

\bigskip

As a consequence of Theorem \ref{Hconv} and the continuous mapping Theorem [see e.g. Theorem 1.3.6 in \cite{vandwell1996}], we obtain the weak convergence of a corresponding Cram\'{e}r-von-Mises and Kolmogorov-Smirnov type test statistic, i.e.
\begin{align}
	\Tb_{n,\tL_2} =\int_{[0,1]^3} \left\{ \Hb_n(x,y,z)\right\}^2\, d(x,y,z) \label{T_nL}
		&\weak \Tb_{C,\tL_2} =\int_{[0,1]^3} \left\{ \Hb_C(x,y,z)\right\}^2\, d(x,y,z), \\
	\Tb_{n,\tKS} = \sup_{[0,1]^3} \left| \Hb_n(x,y,z) \right|   \label{T_nK}
		&\weak \Tb_{C,\tKS}= \sup_{[0,1]^3} \left|  \Hb_C(x,y,z)\right|,
\end{align}
which can be used to construct an asymptotic test for the hypothesis of associativity. Since $\Tb_{n,\tM} \Pconv\infty$ [$\tM\in\{\tL_2,\tKS\}$] if the copula is not associative the null hypothesis  should be rejected for unlikely large values of $\Tb_{n,\tM}$. This gives rise to the demand for critical values of $\Tb_{C,\tM}$ which can be obtained by multiplier bootstrap methods as described in the subsequent paragraph.

\subsection{A multiplier bootstrap approximation} \label{subsec:mult}

It is the purpose of this Section to provide a bootstrap approximation for the distribution
of the limiting variables $\Tb_{C,\tM}$ whose variances depend on the unknown copula in a complicated manner. We begin with an approximation of the distribution of the limiting process
$\Hb_C$. For this purpose we rewrite  the decomposition of the process $\Gb_C$ defined
in \eqref{G_C} as
\begin{align} \label{H_C}
	\Hb_C(x,y,z) &= \Bb_C(x,C(y,z)) - \dceins(x,C(y,z))\Bb_C(x,1) - \dczwei(x,C(y,z)) \Bb_C(1,C(y,z)) \nonumber\\
	& \quad - \left\{ \Bb_C(C(x,y),z) - \dceins(C(x,y),z) \Bb_C(C(x,y),1) - \dczwei(C(x,y),z) \Bb_C(1,z) \right\} \nonumber \\
	& \quad + \dczwei(x,C(y,z)) \left\{ \Bb_C(y,z) - \dceins(y,z)\Bb_C(y,1) - \dczwei(y,z) \Bb_C(1,z) \right\} \nonumber \\
	& \quad + \dceins(C(x,y),z) \left\{ \Bb_C(x,y) - \dceins(x,y)\Bb_C(x,1) - \dczwei(x,y) \Bb_C(1,y) \right\}.
\end{align}
In the following discussion the symbol
\begin{equation}\label{weakcond}
G_n\weakP G
\end{equation}
 denotes weak convergence
in some metric space $\Db$ conditionally on the data in probability
[see  \cite{kosorok2008}]. More precisely, \eqref{weakcond}
holds for random variables $G_n=G_n(\vect{X}_1,\dots,\vect{X}_n,\xi_1,\dots\xi_n),$ $G\in\Db$ if and only if
\begin{align}\label{BL}
    \sup_{h \in BL_1(\Db)} | \Exp_\xi h(G_n) - \Exp h(G)| \Pconv 0
\end{align}
and
\begin{align}\label{am}
    \Exp_\xi h(G_n)^*-\Exp_\xi h(G_n)_* \Pconv 0 \quad \text{for every } h \in BL_1(\Db),
\end{align}
where
\begin{align*}
    BL_1(\Db) =\left\{ f: \Db \rightarrow \R\; | \; ||f||_\infty \leq 1, |f(\beta)-f(\gamma)| \leq d(\beta,\gamma)\;\forall\;\gamma,\beta\in \Db \right\}
\end{align*}
denotes the set of all Lipschitz-continuous functions bounded by $1$.
The subscript $\xi$ in the expectations in \eqref{BL} and   \eqref{am} indicates the
 conditional expectation
with respect to the weights $\xi=(\xi_1,\dots,\xi_n)$ given the data and $h(G_n)^*$
and $ h(G_n)_*$ denote measurable majorants and minorants with respect
to the joint data, including the weights $\xi$. Note also that condition (\ref{BL}) is
motivated by the metrization of weak convergence by the bounded
Lipschitz-metric, see e.g. Theorem 1.12.4 in \cite{vandwell1996}.

The process $\Bb_C$ can be approximated by multiplier bootstrap methods, see \cite{buecher2011, buecdett2010, remiscai2009, segers2011}.
More precisely, let $\xi_1,\dots\xi_n$ denote independent identically distributed random variables with mean $0$ and variance $1$ such that
\begin{equation}\label{xiass}
	||\xi_i||_{2,1}=\int_0^\infty \sqrt{\Prob(|\xi_i|>x)} \, dx < \infty,
\end{equation}  and consider the process
\begin{equation}\label{alphaxi}
\alpha_n^\xi= \sqrt{n}(\ec^{\xi} - \ec),
\end{equation}
 where
\begin{align*}
	\ec^\xi(\vect{x}) = n^{-1}\sum_{i=1}^n \frac{\xi_i}{\bar \xi_n} \ind\{ X_{i1} \leq F_{n1}^-(x_1), X_{i2} \leq F_{n2}^-(x_2) \}
\end{align*}
denotes a multiplier bootstrap version of the estimator.
It was shown in \cite{buecdett2010} and in more detail in \cite{buecher2011} that
$$
\alpha_n^\xi \weakP \Bb_C
$$
i.e. the process $ \alpha_n^\xi$ defined in \eqref{alphaxi}  {converges weakly
to $\Bb_C$ in $l^\infty([0,1]^2)$ conditionally on the data in probability} in the sense of \cite{kosorok2008}.

For the approximation of the partial derivatives in \eqref{H_C} let $\dchatp$ be some estimator of $\dcp$; for instance
an estimator based on the differential quotient as in \cite{remiscai2009} defined by
%For the approximation of the partial derivatives in \eqref{H_C} we utilise the following simple estimators
%proposed by \cite{remiscai2009} and more precisely by \cite{segers2011}
\begin{align}	\label{derivest1}
\dchateins(\vect{u}) &:=
	\begin{cases}
		\frac{\ec(u_1+h,u_2)-\ec(u_1-h,u_2)}{2h} & \text{ if } u_1\in [h,1-h] \\
		\frac{\ec(2h,u_2)}{2h} & \text{ if } u_1 \in[0,h) \\
		\frac{u_2 - \ec(1-2h,u_2)}{2h} & \text{ if } u_1 \in(1-h,1]
	\end{cases}	\\
\dchatzwei(\vect{u}) &:=
	\begin{cases}
		\frac{\ec(u_1,u_2+h)-\ec(u_1,u_2-h)}{2h} & \text{ if } u_2\in [h,1-h] \\
		\frac{\ec(u_1,2h)}{2h} & \text{ if } u_2 \in[0,h) \\
		\frac{u_1 - \ec(u_1,1-2h)}{2h} & \text{ if } u_2 \in(1-h,1]
	\end{cases}   \label{derivest2}
\end{align}
where $h =h_n \rightarrow0$ such that $\inf_n h_n\sqrt{n}>0$ [for a smooth version of these estimators see \cite{scaillet2005}].

\begin{theorem} \label{theo:boot}
Assume that there exists a constant $K$ such that $\|\dchatp\|_\infty\leq K$ for all $n\in\N, p=1,2$ and that
\begin{align*}
	\sup_{\vect{u}\in[0,1]^2:u_p\in[\delta,1-\delta]} \left| \dchatp(\vect{u})-\dcp(\vect{u})\right| \Pconv 0
\end{align*}
for all $\delta\in(0,1/2)$.
If moreover the copula $C$ satisfies Condition \ref{cond:pd} and if the multipliers $\xi_i$ satisfy \eqref{xiass}, then the multiplier process $\Hb_n^\xi$ defined as
\begin{align*}
	\Hb_n^\xi(x,y,z) &= \alpha_n^\xi(x,\ec(y,z)) - \dchateins(x,\ec(y,z))\alpha_n^\xi(x_1,1) - \dchatzwei(x,\ec(y,z)) \alpha_n^\xi(1,\ec(y,z)) \\
	& \quad - \left\{ \alpha_n^\xi(\ec(x,y),z) - \dchateins(\ec(x,y),z)\alpha_n^\xi(\ec(x,y),1) - \dchatzwei(\ec(x,y),z) \alpha_n^\xi(1,z) \right\} \\
	& \quad + \dchatzwei(x,\ec(y,z)) \left\{ \alpha_n^\xi(y,z) - \dchateins(y,z)\alpha_n^\xi(y,1) - \dchatzwei(y,z) \alpha_n^\xi(1,z) \right\} \\
	& \quad + \dchateins(\ec(x,y),z) \left\{ \alpha_n^\xi(x,y) - \dchateins(x,y)\alpha_n^\xi(x,1) - \dchatzwei(x,y)\alpha_n^\xi(1,y) \right\}
\end{align*}
 converges weakly to the process $\Hb_C$ conditional on the data in probability, i.e. $\Hb_n^\xi \weakP \Hb_C$.
\end{theorem}

\textbf{Proof.} Define the process $\tilde \Hb_n^\xi$ by substituting the estimators $\dchateins, \dchatzwei$ and $\ec$ in the definition of $\Hb_n^\xi$ by the true but unknown objects $\dceins, \dczwei$ and $C$. By Lemma A.1 in \cite{buecher2011} it suffices to show that
\begin{align*}
	\|  \Hb_n^\xi-  \tilde \Hb_n^\xi \| _\infty = \sup_{(x,y,z)\in[0,1]^3} | \Hb_n^\xi(x,y,z) -  \tilde \Hb_n^\xi(x,y,z)| \Pconv 0.
\end{align*}
Using the triangle inequality we have to estimate the following 12 summands
\begin{align*}
	\|   \Hb_n^\xi-  \tilde \Hb_n^\xi \| _\infty
	&\leq  \|   \alpha_n^\xi(x,\ec(y,z)) - \alpha_n^\xi(x,C(y,z)) \|_\infty \\
		&+ \|   \dchateins(x,\ec(y,z))\alpha_n^\xi(x,1) - \dceins(x,C(y,z))\alpha_n^\xi(x,1) \|_\infty \\
		&+ \| \dchatzwei(x, \ec(y,z)) \alpha_n^\xi(1,\ec(y,z)) - \dczwei(x, C (y,z))\alpha_n^\xi(1,C(y,z)) \|_\infty \\
		&+ \|  \alpha_n^\xi(\ec(x,y),z) -  \alpha_n^\xi(C(x,y),z) \|_\infty \\
		&+ \|  \dchateins(\ec(x,y),z)\alpha_n^\xi(\ec(x,y),1) - \dceins(C(x,y),z)\alpha_n^\xi(C(x,y),1) \|_\infty \\
		&+ \|   \dchatzwei(\ec(x,y),z) \alpha_n^\xi(1,z) -  \dczwei(C(x,y),z) \alpha_n^\xi(1,z) \|_\infty \\
		&+ \|  \dchatzwei(x,\ec(y,z)) \alpha_n^\xi(y,z) - \dczwei(x,C(y,z)) \alpha_n^\xi(y,z) \|_\infty \\
		&+ \|  \dchatzwei(x,\ec(y,z)) \dchateins(y,z)\alpha_n^\xi(y,1) - \dczwei(x,C(y,z)) \dceins(y,z)\alpha_n^\xi(y,1) \|_\infty \\
		&+ \|  \dchatzwei(x,\ec(y,z))  \dchatzwei(y,z) \alpha_n^\xi(1,z) - \dczwei(x,C(y,z))  \dczwei(y,z) \alpha_n^\xi(1,z) \|_\infty \\
		&+ \|   \dchateins(\ec(x,y),z)  \alpha_n^\xi(x,y) -  \dceins(C(x,y),z)  \alpha_n^\xi(x,y) \|_\infty \\
		&+ \|   \dchateins(\ec(x,y),z) \dchateins(x,y)\alpha_n^\xi(x,1) -  \dceins(C(x,y),z) \dceins(x,y)\alpha_n^\xi(x,1) \|_\infty \\
		&+ \|   \dchateins(\ec(x,y),z)\dchatzwei(x,y)\alpha_n^\xi(1,y)  -  \dceins(C(x,y),z) \dczwei(x,y)\alpha_n^\xi(1,y) \|_\infty,
\end{align*}
 of which one of the hardest cases will be considered exemplarily in the following, namely the third summand
\begin{align*}
	\sup_{(x,y,z)\in[0,1]^3} \left| \dchatzwei(x, \ec(y,z)) \alpha_n^\xi(1,\ec(y,z)) - \dczwei(x, C (y,z))\alpha_n^\xi(1,C(y,z)) \right|.
\end{align*}
 The treatment of the other summands is similar and is omitted for the sake of brevity. We estimate
\begin{align*}
	&\left|\dchatzwei(x, \ec(y,z)) \alpha_n^\xi(1,\ec(y,z)) - \dczwei(x, C (y,z))\alpha_n^\xi(1,C(y,z))\right| \\
	\leq &\left| \dchatzwei(x, \ec(y,z)) - \dczwei(x, \ec(y,z)) \right| \times \left|\alpha_n^\xi(1,\ec(y,z))\right|  \\
	&\quad+ \left|  \dczwei(x, \ec(y,z)) - \dczwei(x, C(y,z)) \right| \times \left| \alpha_n^\xi(1,\ec(y,z))\right|\\
	&\quad+ \left|\dczwei(x,C(y,z))\right| \times\left|\alpha_n^\xi(1,\ec(y,z)) -\alpha_n^\xi(1,C(y,z))  \right| \\
	=:& A_1(x,y,z) + A_2(x,y,z) + A_3(x,y,z)
 \end{align*}
and consider each term separately. For arbitrary $\varepsilon>0$ and $\delta\in(0,1/2)$ we estimate
\begin{align*}
%	\Prob\left( \sup_{\vect{x}\in[0,1]^3} A_1(\vect{x}) \right) \leq \Prob\left( \sup_{\vect{x}\in[0,1]^3, \ec(x_2,x_3)\in[\delta, 1-\delta]} A_1(\vect{x}) \right) + \Prob\left( \sup_{\vect{x}\in[0,1]^3,\ec(x_2,x_3)\notin[\delta, 1-\delta]} A_1(\vect{x} \right)
	\Prob\left( \sup A_1(x,y,z)>\varepsilon \right) \leq \Prob\left( \sup_{\ec(y,z)\in[\delta, 1-\delta]} A_1(x,y,z)>\varepsilon \right) + \Prob\left( \sup_{\ec(y,z)\notin[\delta, 1-\delta]} A_1(x,y,z)>\varepsilon \right)
\end{align*}
where we suppressed the index $(x,y,z)\in[0,1]^3$ at the suprema.
The first probability can be made arbitrary small by the assumptions on $\dchatzwei$ and by the asymptotic tightness of
the process  $\alpha_n^\xi$, see Theorem 2.3 in \cite{buecher2011}.
For the second summand use uniform boundedness of $\dchatzwei$ and the fact that the (unconditional) limit process $\Bb_C(1,\cdot)$ of $\alpha_n^\xi(1,\cdot)$ is a standard Brownian bridge having continuous trajectories which vanish
at $0$ and $1$. By decreasing $\delta$ the probability can be made
 arbitrary small, see \cite{segers2011} for an rigorous treatment of this argument.

Since $\dczwei$ is uniformly continuous if the second coordinate is bounded away from zero and one the second summand $A_2(x,y,z)$ can be treated similarly. Regarding $A_3(x,y,z)$ note that $\alpha_n^\xi$ is asymptotically uniformly equicontinuous
[Theorem 2.3 in \cite{buecher2011}]
and that $\sup_{(y,z)\in[0,1]^2}|\ec(y,z)-C(y,z)|\Pconv0$ which yields
$$\sup_{(y,z)\in[0,1]^2}\left|\alpha_n^\xi(1,\ec(y,z)) -\alpha_n^\xi(1,C(y,z))  \right|\Pconv0.$$
By boundedness of $\dczwei$ this yields the assertion  $\sup_{(x,y,z)\in[0,1]^3}A_3(x,y,z)\Pconv0$. \qed

\bigskip

\begin{remark}\label{remark} ~~\\
 (a) Note that the assumptions on the estimator $\dchatp$ for the partial derivatives $\dcp$ are e.g. satisfied for the estimators defined in \eqref{derivest1} and \eqref{derivest2}, see Lemma 4.1 in \cite{segers2010}.

\medskip

(b) 
Note that Theorem \ref{theo:boot}  holds independently of the hypothesis of associativity.
As a consequence of the continuous mapping theorem for the bootstrap, see Proposition 10.7 in \cite{kosorok2008}, we can conclude that
\begin{align} \label{tnb}
	\Tb_{n,\tL_2}^\xi=\int_{[0,1]^3} \left\{ \Hb_n^\xi(x,y,z)\right\}^2\, d(x,y,z) \weakP \Tb_{C,\tL_2}, \quad \Tb_{n,\tKS}^\xi = \sup _{[0,1]^3} \left| \Hb_n^\xi(x,y,z)\right| \weakP \Tb_{C,\tKS}
\end{align}
and the latter convergence suggests to use the following approach in order to obtain an asymptotic level-$\alpha$ test for the hypothesis of associativity.
\begin{enumerate}
	\item Compute the statistic $\Tb_{n,\tM}$ [$\tM\in\{\tL_2,\tKS\}$].
	\item Choose the  number of bootstrap replications $B\in\N$. For $b=1,\ldots ,B$
	 simulate independent replications of the random variables $\xi_1,\dots,\xi_n$
	 and denote the result form the $b$-th iteration by $\xi_{1,b},\dots,\xi_{n,b}$.
	\item For $b=1,\ldots, B$ compute  the statistics $\Tb_{n,\tM}^{(\xi,b)}$ defined in \eqref{tnb}  from the data
$\vect{X}_1,\dots,\vect{X}_n$ and  the multipliers $\xi_{1,b},\dots,\xi_{n,b}$
and determine the ${(1-\alpha)}$-quantile $q_{1-\alpha,\tM}^\xi$ of the empirical distribution of the sample $\{\Tb_{n,\tM}^{(\xi,b)}\}_{b=1,\dots,B}$.
\item Reject the null hypothesis of associativity whenever $\Tb_{n,\tM}>q_{1-\alpha,\tM}^\xi$
\end{enumerate}
Since $\Tb_{n,\tM}\Pconv\infty$ and $\Tb_{n,\tM}^\xi\weakP\Tb_{C,\tM}$ if the copula is not associative the test is consistent against all alternatives satisfying the conditions of Theorem \ref{theo:boot}.
\end{remark}

\section{Testing Archimedeanity}\label{sec:testarch}
\def\theequation{3.\arabic{equation}}
\setcounter{equation}{0}

\begin{figure}
    \includegraphics[width=0.95\textwidth]{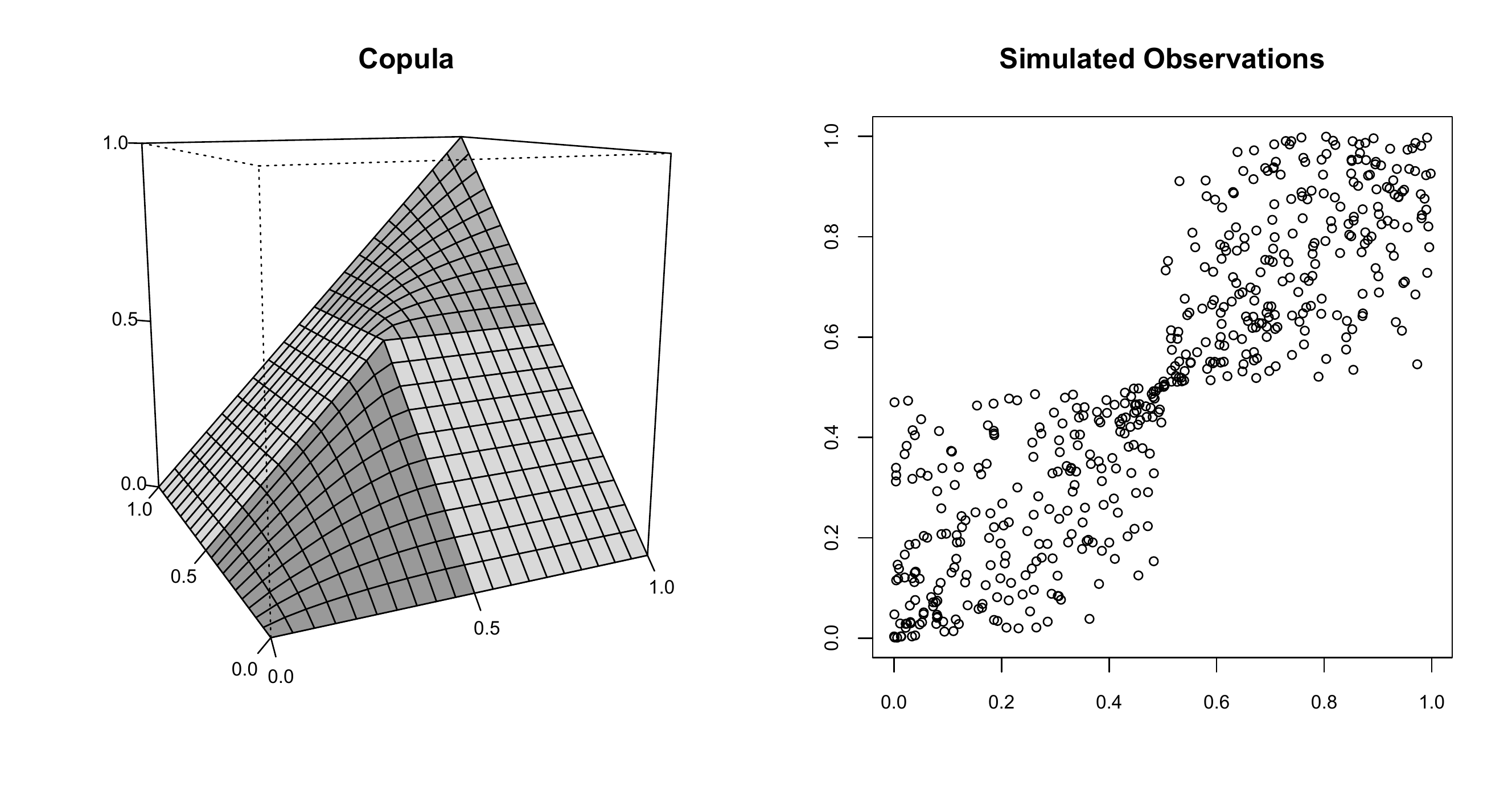}
  \caption{\label{pic:os}
{ \textit{Left picture: Ordinal sum copula. Right picture: $500$ corresponding simulated observations.  }
} }
\end{figure}

As stated in the Introduction a bivariate copula $C$ is Archimedean if and only if $C$ is an associative copula satisfying $C(u,u)<u$ for all $u\in(0,1)$. Associativity has been dealt with in the preceding paragraph and it remains to handle non-Archimedean copulas which my be associative but satisfy $C(q,q)=q$ for some $q\in(0,1)$. Due to Theorem 1 in \cite{jaworski2010} or by the results in Section 2.4 of \cite{alsfrasch2006} all those copulas may be expressed as an ordinal sum of Archimedean copulas. An ordinal sum copula is defined as following [cf. Section 3.2.2 in \cite{nelsen2006}]: let $\{J_i\}_{i\in I}$ be a countable partition of non-overlapping closed intervals $J_i=[a_i,b_i]$ whose union is $[0,1]$. If moreover $\{C_i\}_{i\in I}$ is a collection of copulas, then the ordinal sum of $\{C_i\}_{i \in I}$ with respect to $\{J_i\}_{i \in I}$ is the copula $C$ defined by
\begin{align*}
	C(\vect{u}) = \begin{cases} a_i+(b_i-a_i)C_i\left(\frac{u_1-a_i}{b_i-a_i}, \frac{u_2-a_i}{b_i-a_i}\right) & \text{if } \vect{u}\in J_i\times J_i\\
		\min\{u_1,u_2\} & \text{otherwise.}
	\end{cases}
\end{align*}
Note that $C(b_i,b_i)=b_i$ for all $b_i$ and that ordinal sum copulas put no mass on $[0,1]^2\setminus \bigcup_{i\in I} J_i\times J_i$. In Figure \ref{pic:os} we illustrate the ordinal sum   of a Gumbel copula $C_1$ with parameter $\theta_1=1.5$ and a Clayton Copula $C_2$ with parameter $\theta_2=1$, where $J_1=[0,1/2], J_2=[1/2,1]$. Note that Kendall's $\tau$ of $C$ is equal to $2/3$, while it equals $1/3$ for both $C_1$ and $C_2$.

In order to check for $C(q,q)=q$ for some $q\in(0,1)$ we propose the following modification of the statistic $\Tb_{n,\tM}$
\begin{align*}
	\Sb_{n,\tM} = \Tb_{n,\tM} + k_n \phi(A_n(\ec)),
\end{align*}
where   $\Tb_{n,\tM}$ is defined in \eqref{T_nL} and \eqref{T_nK},
$k_n\sim n^\alpha, \alpha\in(0,1/2)$ is some constant chosen by the statistician, $\phi$ is some increasing function with $\phi(0)=0$ and
\begin{align*}
	A_n(\ec) =\max\left\{ \frac{i}{n} \left( 1-\frac{i}{n} \right) \, : \, \ec \left( \frac{i}{n},\frac{i}{n} \right) = \frac{i}{n} \right\}.
\end{align*}
Intuitively, $A_n(\ec)$ should be ``large'' for copulas which satisfy $C(q,q)=q$ for some $q\in(0,1)$. For a decent choice of $k_n$ and $\phi$ we refer the reader to Section \ref{sec:sim}.

\begin{figure}[ht]
    \includegraphics[width=0.95\textwidth]{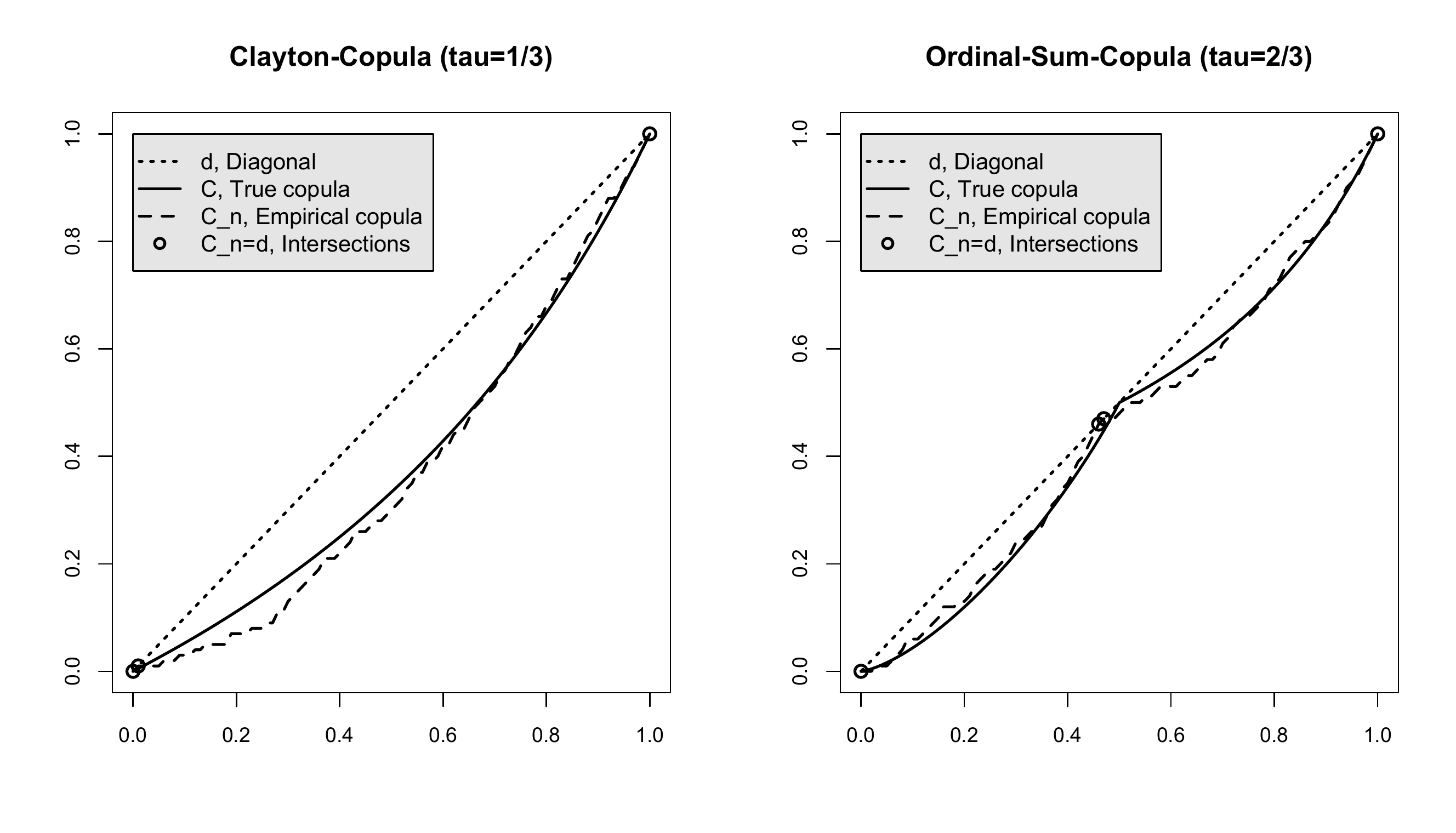}
  \caption{\label{pic:diag}
\textit{The solid lines show the diagonal section of a Clayton copula (left) and an ordinal sum copula (right), while the dashed line show one realization of the corresponding empirical copula. The circled points mark the locations where $\ec(i/n,i/n)=i/n$.}
} 
\end{figure}

In Figure \ref{pic:diag}
we illustrate the points $i/n$ at which $\ec(i/n,i/n)=i/n$ for two specific examples, the Clayton copula with $\theta=1$ and the ordinal sum copula depicted in Figure \ref{pic:os}. The solid and dashed lines correspond to the true copula and the empirical copula calculated for a set of $n=100$ observations, respectively. For the ordinal sum copula there  always exist some points $i/n$ in a neighbourhood of $1/2$ such that $\ec(i/n,i/n)=i/n$, see the proof of the following Proposition, which is sufficient for the derivation of the asymptotic properties of the statistic
$\Sb_{n,\tM}$.

\begin{prop} \label{prop:an} ~
\begin{enumerate}[a)]
\item Suppose $C$ is an Archimedean copula satisfying Condition \ref{cond:pd} and that the coefficients of tail dependence
	\begin{align*}
		\lambda_L = \lim_{u\rightarrow 0} \frac{C(u,u)}{u} \quad \text{ and } \quad \lambda_U= \lim_{u\to 1} \frac{1-2u+C(u,u)}{1-u}
	\end{align*}
	exist and are smaller than $1$. Then it holds
	$$A_n(\ec) = o_\Prob(n^{-\alpha})$$
	for any $\alpha<1/2$.
\item If there exists a $q\in(0,1)$ such that $C(q,q)=q$ then it holds
	$$A_n(\ec)\geq q(1-q) + o_\Prob(1).$$
\end{enumerate}
\end{prop}

\textbf{Proof. } We start with the proof of \textit{a)}. First choose $\delta>0$ and $\lambda<1$ such that
$$\frac{C(u,u)}{u}\vee \frac{1-2u+C(u,u)}{1-u} \leq\lambda  $$
for all $u\in[0,\delta]\cup [1-\delta,1]$ and use the decomposition
\begin{align}\label{eq:deca}
	A_n(\ec) = A_n(\ec, [0,\delta)) + A_n(\ec,[\delta, 1-\delta]) + A_n(\ec,(1-\delta,1]),
\end{align}
where
$$
A_n(\ec,B) = \max\left\{ i/n \left( 1-i/n \right) \, : \, \ec \left( i/n, i/n \right) = i/n \text{ and } i/n \in B \right\}
$$
 for some set $B \subset [0,1]$  (with the convention that $\max\emptyset=0$).
 Consider each term separately and define $M_n=\sup_{u\in[0,1]} \left| \ec(u,u) - C(u,u)\right|,$ which is of order $O_\Prob(n^{-1/2})$ under Condition \ref{cond:pd}. Now let $i/n\in(0,\delta)$ be such that $\ec(i/n,i/n)=i/n$.  Due to the estimate
\begin{align*}
	i/n(1-\lambda) \leq i/n\left( 1- \frac{C(i/n,i/n)}{i/n} \right) = i/n - C(i/n,i/n) = \ec(i/n,i/n)-C(i/n,i/n) \leq M_n
\end{align*}
we have $i/n(1-i/n) \leq i/n \leq M_n/(1-\lambda)$ and we can conclude that
\begin{equation}\label{a1}
A_n(\ec,[0,\delta)) \leq \frac{M_n}{1-\lambda} =O_\Prob(n^{-1/2}).
\end{equation}
A similar calculation shows that for $i/n\in(1-\delta,1]$ with $\ec(i/n,i/n)=i/n$ we have $(1-i/n)(1-\lambda) \leq M_n$ which in turn implies
\begin{equation}\label{a2}
A_n(\ec,(1-\delta,1]) \leq \frac{M_n}{1-\lambda} =O_\Prob(n^{-1/2}).
\end{equation}
It remains to estimate the second summand $A_n(\ec,[\delta, 1-\delta])$ of the decomposition \eqref{eq:deca}. For continuity reasons we can choose a $\kappa>0$ such that $u-C(u,u) \geq \kappa$ for all $u\in[\delta,1-\delta]$. If there was a $q\in[\delta,1-\delta]$ such that $\ec(q,q)=q$, it would follow that $M_n\geq \ec(q,q)-C(q,q) \geq \kappa$ and therefore we have for any $\varepsilon>0$
$$\Prob\left(n^\alpha A_n(\ec,[\delta, 1-\delta]\right) > \varepsilon) \leq \Prob( \exists \, q\in[\delta,1-\delta]: \ec(q,q)=q) \leq \Prob( M_n\geq \kappa) \Pconv 0.$$
A combination of \eqref{a1} and \eqref{a2} with this result proves part  \textit{a)} of the proposition.

For the proof of part \textit{b)} let $n_1=\# \{1\leq i \leq n :  (F_{1}(X_{i1}), F_{2}(X_{i2})) \in [0,q]^2 \}$. Since $C(q,q)=q$ implies that the mass of $C$ is concentrated on $(0,q)^2 \cup (q,1)^2$ we have $(F_1(X_{i1}), F_2(X_{i2})) \in [0,q]^2$ if and only if  $X_{i1}\leq X_{n_1:n,1}$ and  $X_{i2}\leq X_{n_1:n,2}$,
%\begin{align*}
%%	n_1= \# \{1\leq i \leq n :  X_{i1}\leq X_{n_1:n,1}, X_{i2}\leq X_{n_1:n,2}\},
%\{1\leq i \leq n :  (F_{1}(X_{i1}), F_{2}(X_{i2})) \in [0,q]^2 \}=  \{1\leq i \leq n :  X_{i1}\leq X_{n_1:n,1}, X_{i2}\leq X_{n_1:n,2}\},
%\end{align*}
where $X_{j:n,p}=F_{np}^-(j/n)$ denotes the $j$-th order statistic of $X_{1p},\dots,X_{np}$ (for $p=1,2$). This yields $\ec(n_1/n,n_1/n) = n_1/n$, which entails the assertion by
$$A_n(\ec) \geq \frac{n_1}{n}\left( 1- \frac{n_1}{n} \right) \Pconv q(1-q) > 0.   $$ \qed

\begin{remark}~

\begin{enumerate}[a)]
	\item The conditions on the coefficients of tail dependence in part \textit{a)}
 of Proposition \ref{prop:an}
can be equivalently expressed by conditions on the regular variation of the Archimedean generator of $C$. For a thorough discussion of these issues the reader is referred to the work of \cite{charsege2009}.
	\item Exploiting Theorem G.1 in \cite{genesege2009} and Proposition 4.2 in \cite{segers2011} one can improve the rate of convergence in
	part $\textit{a)}$ of Proposition \ref{prop:an}
	to any $\alpha<3/4$. It is our conjecture that the term is in fact of order $O_\Prob(1/n)$, but we were not able to derive the asymptotic distribution of $nA_n(\ec)$. Since we do not need a refined rate for our purposes here we omit a deeper discussion and defer these issues to future research.
\end{enumerate}
\end{remark}

From now on, suppose that the conditions of Theorem \ref{theo:boot} and  Proposition \ref{prop:an} hold.
We can conclude that $\Sb_{n,\tM}$ weakly converges to $\Tb_{C,\tM}$ if the copula $C$ is Archimedean, while $\Sb_{n,\tM}$ converges to $+\infty$ in probability if $C$ is non-Archimedean, i.e. if it is either non-associative (by the results of Section \ref{sec:test}) or if there exists a $q\in(0,1)$ such that $C(q,q)=q$ (by Proposition \ref{prop:an}). The quantiles of $\Tb_{C,\tM}$ can be approximated by the multiplier method described in Section \ref{subsec:mult}. Analogously to the discussion at the end of the Section \ref{subsec:mult} we can use the multiplier bootstrap to
obtain an asymptotic  level-$\alpha$ test for the hypothesis of Archimedeanity, which
is consistent  against all alternatives satisfying the Condition \ref{cond:pd}.  Its finite sample properties will be investigated in the following section.

\section{Finite sample properties} \label{sec:sim}

We conclude this paper with a  simulation study regarding the finite sample performance of the proposed tests for Archimedeanity and Associativity.
We consider the following six copula models: 
%the Archimedean Gumbel and Clayton copula models, the non-associative $t$-copula and the asymmetric negative logistic copula model [with fixed parameters $\psi_1=2/3, \psi_2=1$, see \cite{joe1990}] and two associative ordinal sum models.
%The latter two models are both based on the partition $J_1=[0,1/2], J_2=[1/2,1]$. The first ordinal sum model will be denoted by $\text{Ordinal}_A$  and is based on the Gumbel  $(C_1)$ and Clayton $(C_2)$ copula, while the second model (denoted by $\text{Ordinal}_B$) is based on two Clayton copulas $(C_1=C_2)$.

\begin{itemize}
\item The Gumbel copula, which is Archimedean. 
\item The Clayton copula, which is Archimedean.
\item The $t$-copula with fixed degree of freedom $df=1$, which is non-associative.
\item The asymmetric negative logistic model [see \cite{joe1990}] with fixed parameters $\psi_1=2/3, \psi_2=1$, which is non-associative.
\item An ordinal sum model based on the partition $J_1=[0,1/2], J_2=[1/2,1]$ and the Gumbel  $(C_1)$ and Clayton $(C_2)$ copula, denoted by $\text{Ordinal}_A$. The model is associative.
\item An ordinal sum model based on the partition $J_1=[0,1/2], J_2=[1/2,1]$ and the two Clayton $(C_1=C_2)$ copulas, denoted by $\text{Ordinal}_B$. The model is associative.
\end{itemize}
The parameters of the models are chosen in such a way that the coefficient of upper tail dependence $\lambda_U$ is either $1/4$ or $1/2$ [for the asymmetric negative logistic model] or that Kendall's-$\tau$ is either $1/3$ or $2/3$ [for the remaining five models]. For $\tau_{\text{Ordinal}_A}=1/3$ (or $2/3$, resp.) we chose $\tau_{\text{Gumbel}}=0$ ($1/3$) and $\tau_{\text{Clayton}}=-2/3$ ($1/3$), while $\tau_{\text{Clayton}}=-1/3$ ($1/3$) for $\tau_{\text{Ordinal}_B}=1/3$ ($2/3$).

We generated $1000$ random samples of sample sizes $n=200$ and $n=500$ and calculated the empirical probability of rejecting the null hypotheses of Archimedeanity or Associativity for $\tM\in\{\tL^2,\tKS\}$. For each sample of size $n=200, 500$ we carried out $B=200$ Bootstrap replications based on the multiplier method, where we chose a $\mathcal{U}(\{0,2\})$-distribution for the multipliers [i.e. $\Prob(\xi=0)=\Prob(\xi=2)=1/2$, s.t. $\mu=\tau=1$] and used $h=n^{-1/4}$ to estimate the partial derivatives. The critical values of the tests are determined by the method described in Section \ref{sec:test}.  The penalty term $\Sb_{n,\tM} - \Tb_{n,\tM} = k_n\phi(A_n(C_n))$ is chosen in the  following data-adaptive way: first of all, we set $\phi(x)=(4x)^2$ in order to give more emphasis to values around the maximal value of $A_n(C_n)$ [which equals $1/4$]. The constant $k_n$ is chosen according to the distribution of the bootstrap approximation: if $q^\xi_{0.05,\tM}$ denotes the $0.05$-quantile of the sample $\{\Tb_{n,\tM}^{\xi,b}\}_{b=1,\dots,B}$ we set $k_n=q^\xi_{0.05,\tM} n^{1/4}$. The latter choice guarantees that under $H_0$ the error term is small compared to the distribution of $\Tb_{C,\tM}$.

The results are stated in Table \ref{tab:test01}. The entries of the table represent the empirical probabilities of rejecting the null hypothesis of Archimedeanity and of Associativity [in brackets] for both the $\tL^2$-test [first two columns] and the $\tKS$-test [last two columns].
We observe that the nominal level of the four tests are accurately approximated for the four Archimedean copulas under investigation. The $\tL^2$-test tends to be more conservative than the $\tKS$-test. Also note that the values for $\Sb_{n,\tM}$ and $\Tb_{n,\tM}$ differ only by a very small amount meaning that the penalty term $k_n\phi(A_n(C_n))$ is of negligible magnitude under the null hypothesis.

The $t$-copula and the asymmetric negative logistic models are non-associative and the results in Table \ref{tab:test01} reveal that these deviations are detected by both tests for Associativity, with better results for the $t$-copula and for stronger dependence [measured by either $\tau$ or $\lambda_U$]. The power properties of the $\tL^2$-test outclass the properties of the $\tKS$-test for all four non-associative models under investigation, such that the former test seems to be generally preferable.

Regarding the (associative) ordinal sum models both tests for associativity are very conservative. Note that the asymptotic theory developed in Section \ref{sec:test} does not apply for these models since the partial derivatives of the corresponding copulas are not continuous on $\{1/2\}\times[0,1]$ and $[0,1]\times\{1/2\}$. Regarding the power properties the $\tKS$-test for Archimedeanity performs slightly better for the ordinal sum alternatives.

\begin{table}
\begin{center}
\begin{tabular}{ l@{ \qquad\quad }   r r c r r }\hline\hline
&\multicolumn{2}{c}{$\tL2$-Test} & & \multicolumn{2}{c}{$\tKS$-Test} \\
%& $\tL2$-Test & & $\tKS$-Test &  \\
 & \multicolumn{1}{c}{0.1} & \multicolumn{1}{c}{0.05}&  \qquad\qquad  & \multicolumn{1}{c}{0.1}  & \multicolumn{1}{c}{0.05} \\ \hline\hline
\multicolumn{1}{l}{\bf n=200:} & &&&&\\ \hline
 Clayton($\tau=1/3$) &    0.071 (0.071) & 0.038 (0.037)  & & 0.088 (0.088) & 0.050 (0.050) \\
  Clayton($\tau=2/3$) &    0.030 (0.016) & 0.011 (0.009)  & & 0.124 (0.078) & 0.068 (0.036) \\
 Gumbel($\tau=1/3$) &    0.079 (0.079) & 0.043 (0.043)  & & 0.082 (0.082) & 0.046 (0.045) \\
Gumbel($\tau=2/3$) &    0.034 (0.032) & 0.015 (0.013)  & & 0.108 (0.098) & 0.065 (0.057) \\  \hline
$t(\tau=1/3, df=1)$  &    0.953 (0.953) & 0.886 (0.884) &  & 0.562 (0.558) & 0.380 (0.376) \\
$t(\tau=2/3, df=1)$  &    0.748 (0.726) & 0.592 (0.564)  & & 0.392 (0.355) & 0.258 (0.231) \\
Aneglog($\lambda_U=0.25$)  &   0.112 (0.112) & 0.061 (0.061)  & & 0.105 (0.105) & 0.059 (0.059) \\
Aneglog($\lambda_U=0.5$)  &    0.641 (0.641) & 0.536 (0.536)  & & 0.363 (0.356) & 0.225 (0.222) \\
$\text{Ordinal}_A(\tau=1/3)$ &    0.996 (0.000) & 0.827 (0.000)  & & 1.000 (0.012) & 1.000 (0.005) \\
$\text{Ordinal}_A(\tau=2/3)$ &    1.000 (0.004) & 1.000 (0.001)  & & 1.000 (0.079) & 1.000 (0.041) \\
$\text{Ordinal}_B(\tau=1/3)$ &    1.000 (0.000) & 1.000 (0.000)  & & 1.000 (0.045) & 1.000 (0.021) \\
$\text{Ordinal}_B(\tau=2/3)$ &    1.000 (0.006) & 1.000 (0.001)  & & 1.000 (0.057) & 1.000 (0.030) \\ \hline
\multicolumn{1}{l}{\bf n=500:} & &&&&\\ \hline
 Clayton($\tau=1/3$) &    0.082 (0.082) & 0.051 (0.051)  & & 0.088 (0.088) & 0.036 (0.036) \\
  Clayton($\tau=2/3$) &    0.062 (0.059) & 0.032 (0.027)  & & 0.090 (0.082) & 0.046 (0.039) \\
 Gumbel($\tau=1/3$) &    0.091 (0.091) & 0.045 (0.045)  & & 0.105 (0.105) & 0.050 (0.050) \\
Gumbel($\tau=2/3$) &    0.072 (0.070) & 0.033 (0.032)  & & 0.124 (0.121) & 0.068 (0.066) \\  \hline
$t(\tau=1/3, df=1)$  &    1.000 (1.000) & 1.000 (1.000) &  & 0.954 (0.953) & 0.871 (0.871) \\
$t(\tau=2/3, df=1)$  &    0.998 (0.990) & 0.998 (0.990)  & & 0.818 (0.811) & 0.655 (0.650) \\
Aneglog($\lambda_U=0.25$)  &    0.237 (0.237) & 0.173 (0.173)  & & 0.124 (0.124) & 0.069 (0.069) \\
Aneglog($\lambda_U=0.5$)  &    0.979 (0.979) & 0.947 (0.947)  & & 0.716 (0.716) & 0.584 (0.584) \\
$\text{Ordinal}_A(\tau=1/3)$ &    1.000 (0.000) & 0.980 (0.000)  & & 1.000 (0.022) & 1.000 (0.007) \\
$\text{Ordinal}_A(\tau=2/3)$ &    1.000 (0.021) & 1.000 (0.009)  & & 1.000 (0.093) & 1.000 (0.038) \\
$\text{Ordinal}_B(\tau=1/3)$ &    1.000 (0.000) & 1.000 (0.000)  & & 1.000 (0.005) & 1.000 (0.023) \\
$\text{Ordinal}_B(\tau=2/3)$ &    1.000 (0.013) & 1.000 (0.004)  & & 1.000 (0.082) & 1.000 (0.037) \\  \hline
 \end{tabular}
 \caption{ {\it Simulated rejection probabilities for the multiplier bootstrap-tests for Archimedeanity and for Associativity (in brackets). The sample size is $n=200$ or $n=500$, $B=200$ Bootstrap-replicates and 1000 simulation runs have been performed. The first four lines are Archimedean copula models, the $t$-models are not associative and the ordinal sum-models are associative, but not Archimedean.} }\label{tab:test01}
\end{center}
\end{table}

\bigskip
\medskip

{\bf Acknowledgements} % The authors would like to thank Martina Stein, who typed parts of this manuscript with considerable technical expertise.
This work has been supported   by the Collaborative
Research Center ``Statistical modeling of nonlinear dynamic processes'' (SFB 823, Teilprojekt A1, C1) of the German Research Foundation (DFG).

\bibliography{archcop}

\end{document}